\documentclass[12pt]{article}
\usepackage{amssymb}
\usepackage{amsmath}
\usepackage{amsfonts}
\usepackage{mathrsfs}
\usepackage{amssymb,amsmath}
\begin{document}
\title{\textbf{Normalized Ground State Solutions for Fractional Schr\"{o}dinger Systems with General Nonlinearities}}

\author{Hichem Hajaiej$^{\mathrm{a,}}$, Linjie Song$^{\mathrm{b,c,d}}$\thanks{%
		Song is supported by "Shuimu Tsinghua Scholar Program".}  \\
%EndAName
\\
{\small $^{\mathrm{a}}$ Department of Mathematics, California State University at Los Angeles, Los Angeles, CA 90032, USA}\\
{\small $^{\mathrm{b}}$ Department of Mathematical Sciences, Tsinghua University, Beijing 100084, China} \\
{\small $^{\mathrm{c}}$Institute of Mathematics, AMSS, Chinese Academy of Science,
	Beijing 100190, China}\\
{\small $^{\mathrm{d}}$University of Chinese Academy of Science,
	Beijing 100049, China}
}
\date{}
\maketitle

\begin{abstract}
Normalized ground state solutions (NGSS) of Schr\"{o}dinger equations (SE) have attracted the attention of many research groups during the last decades. This is essentially due to their relevance in many fields in physics and engineering, where the stable and most attractive solutions happen to be the normalized ones. For a single (SE), recent developments lead to the establishment of existence and non-existence results for a wide range of natural nonlinearities in (SE) in the sub-critical, critical and super-critical regimes. However for systems of (SE), there are still many interesting open questions for basic nonlinearities. It certainly requires innovative ideas to shed some light to treat these complex situations. So far, only a very few specific nonlinearities have been addressed. Unlike the single (SE), the corresponding strict sub-additivity inequality is challenging and an improved concentration-compactness theorem is critical to treat (NGSS) of systems of (SE). The aim of this paper is to establish the existence of (NGSS) for a large class of nonlinearities. This class includes many relevant pure-power type nonlinearities and can be easily extended to Hartree type nonlinearities (and a combination of both). The presence of the fractional Laplacian adds a considerable difficulty to rule out the dichotomy. We were able to overcome this challenge and establish very general assumptions ensuring the strict subadditivity of the constrained energy functional. We believe that our approach will open the door to many unresolved problems.
\end{abstract}

{\small Keywords:} {\small  Fractional Laplacian, Vectorial Schr\"{o}dinger, constrained minimization problem}

{\small Mathematics Subject Classification (2010):} 49K99

{\small Data availability statement:} My manuscript has no associate data

%{\small Abbreviated title: }

\newtheorem{theorem}{Theorem}[section]
\newtheorem{remark}[theorem]{Remark}
\newtheorem{proposition}[theorem]{Proposition}
\newtheorem{lemma}[theorem]{Lemma}
\newtheorem{definition}[theorem]{Definition}
\renewcommand{\theequation}{\thesection.\arabic{equation}}
\catcode`@=11 \@addtoreset{equation}{section} \catcode`@=12
\allowdisplaybreaks
\section{Introduction}

In this paper, we consider the following coupled system:
\begin{align}\label{1.1}
(-\triangle)^{\alpha}u_i+\lambda_i u_i+\partial_i F(x,u_1,\ldots,u_m)=0,\ i=1,\ldots,m,
\end{align}
where $x\in\mathbb{R}^N$, $m\geq 2,$ $\alpha>\frac{N}{2}.$

A standard way to obtain some solutions to the above system is to minimize the associated energy functional: (under a suitable constraint)
\begin{align*}
J(\textbf{u})=\frac{1}{2}\int_{\mathbb{R}^N}\sum_{i=1}^m|\Lambda^\alpha u_i|^2dx-\int_{\mathbb{R}^N}F(x,\textbf{u})dx,
\end{align*}
where $\textbf{u}=(u_1,\ldots,u_m)$. $F(x,\textbf{u})=F(x,u_1,\ldots,u_m).$ $\Lambda=(-\triangle)^{\frac{1}{2}}$ is the Zygmund operator and $(-\triangle)^\alpha=\Lambda^{2\alpha}$ can be defined through the Fourier transform:
\begin{align*}
\widehat{\Lambda^{2\alpha} f}(\xi)=|\xi|^{2\alpha} \hat{f}.
\end{align*}The norm of the fractional order Hilbert space $H^{\alpha}(\mathbb{R}^N)$ is
defined by
$$
\|\cdot\|_{H^{\alpha}}=\|\mathcal{F}^{-1}[(1 + |\xi|^{\alpha})^{2}\mathcal{F}(\cdot)]\|_{L^{2}}.
$$
Another definition of the fractional Laplacian is given by the formula:
\begin{align*}
(-\triangle)^\alpha u=C_\alpha P.V.\int_{\mathbb{R}^N}\frac{u(x)-u(y)}{|x-y|^{1+2\alpha}}dy.
\end{align*}
where $u(x):\mathbb{R}^N\mapsto \mathbb{C}$ is smooth enough.

A minimizer of the multi-constrained variational problem \eqref{1.2} below is called a \textbf{normalized ground state solution (NGSS)} of the system \eqref{1.1}
\begin{align}\label{1.2}
I_{c_1,\ldots,c_m}=\inf\{J(\textbf{u}):\textbf{u}\in S_c\},
\end{align}
\begin{align*}
S_c=\{\textbf{u}\in H^\alpha(\mathbb{R}^N)\times\cdots\times H^\alpha(\mathbb{R}^N):\int_{\mathbb{R}^N}|u_i|^2dx=c_i,1\leq i\leq m\},
\end{align*}
where $c_1,\ldots,c_m>0$, $c= (c_1,\ldots,c_m)$. Under some additional regularity assumptions on $F$, a solution of \eqref{1.2} is also a solution of \eqref{1.1}, where $\lambda_i$ are Lagrange multipliers.

%\sum_{i=1}^m c_i$.
For the convenience of the reader, we abbreviate $\int_{\mathbb{R}^N}$ by $\int$ in the sequel and we denote $ H^\alpha(\mathbb{R}^N)\times\cdots\times H^\alpha(\mathbb{R}^N)$ by $H^\alpha$. Moreover, $C$ denotes a generic positive constant which may varies in different estimates. $C_{a,b,\ldots}>0$ is also a generic constant which depends on $a,b,\ldots$.

To the best of our knowledge, there are no previous contributions on the nonlinearity $F$ without standard convexity, compactness, symmetry and monotonicity properties, except for \cite{HS} where the authors studied the system \eqref{1.1} in the classical setting $(\alpha=1)$, and considered nonlinearities $F$ that are products of functions involving all the variables. In this work, we study a more general operator and we will address a wide range of nonlinearities. Hypotheses $(H_3)$ and $(H_6)$ below were the key ideas to include a wide range of functions $F$ that are not necessary the product of functions. To the best of our knowledge, this type of assumptions seem to be novel, and was critical to obtain the sub-additivity condition for general non-linearities. 

 More precisely, our goal is to establish the existence of minimizers of \eqref{1.2} for a given Carath$\mathrm{\acute{e}}$odory function $F$ satisfying the following assumptions:
\newline
($H_1$): For all $x\in \mathbb{R}^N, \textbf{u}\in \mathbb{R}^m$, there exist $A,B>0$ and $0<l_{1,i}<\frac{4\alpha}{N}$ such that
\begin{align*}
&0\leq F(x,u_1,\ldots,u_m)\leq A\Big(\sum_{i=1}^m|u_i|^2+\sum_{i=1}^m|u_i|^{l_{1,i}+2}\Big),\\
&\partial_jF(x,u_1,\ldots,u_m)\leq B\Big(\sum_{i=1}^m|u_i|+\sum_{i=1}^m|u_i|^{l_{1,i}+1}\Big), \ \forall\ 1\leq j\leq m,
\end{align*}
where $\partial_j:=\partial_{u_j}$.
\newline
($H_2$): There exist $\beta>0,$ $S>0,$ $R>0$, $s_1,\ldots,s_m>0$, $s=\sum_{i=1}^ms_i$, $t\in[0,2)$ such that for all $|x|\geq R$ and $\sum_{i=1}^m|u_i|^2<S$ with $t<N(1-\frac{s}{2\alpha})+2$, it holds
$$
F(x,\textbf{u})>\beta|x|^{-t}|u_1|^{s_1}\cdots|u_m|^{s_m}.
$$
($H_3$): For all $j \in\{0,1,\cdots,m\}$, we can write $F(x,u_1,\cdots,u_m) = F_{1j}(x,u_1,\cdots,u_m) + F_{2j}(x,u_1,\cdots,u_{j-1},u_{j+1},\cdots,u_m)$. Furthermore, for all $x\in \mathbb{R}^N,\textbf{u}\in\mathbb{R}^m$ and $\theta_j\geq 1$, we have:
$$
F_{1j}(x,u_1,u_2,\cdots,\theta_ju_j,\cdots,u_m)\geq\theta_j^{2}F_{1j}(x,u_1.u_2,\cdots,\theta_ju_j,\cdots,u_m),
$$
and there exists $L\in\mathbb{R}^N$ such that $$F_{2j}(x+L,u_1,\cdots,u_{j-1},u_{j+1},\cdots,u_m)=F_{2j}(x,u_1,\cdots,u_{j-1},u_{j+1},\cdots,u_m)$$ for all $x\in\mathbb{R}^N$ and $\textbf{u}\in\mathbb{R}^m.$
%For all $x\in \mathbb{R}^N,\textbf{u}\in\mathbb{R}^m$ and $\theta_i\geq 1,i=1,\ldots,m$, we have:
%$$
%F(x,\theta\textbf{u})\geq\theta_{max}^{2}F(x,\textbf{u}),
%$$
%where $\theta_{max}=\max_{1\leq i\leq m}\theta_i$.

In addition, we assume that there exists a periodic function $F^\infty(x,\textbf{u})$; that is there exists $L\in\mathbb{R}^N$ such that $F^\infty(x+L,\textbf{u})=F^\infty(x,\textbf{u})$ for all $x\in\mathbb{R}^N$ and $\textbf{u}\in\mathbb{R}^m,$ satisfying $(H_2)$ and the following assumptions:
\newline
($H_4$): There exists $0<l_{2,i}<\frac{4\alpha}{N}$ such that it holds uniformly for all $\textbf{u}\in\mathbb{R}^m$,
$$
\lim_{|x|\rightarrow+\infty}\frac{F(x,\textbf{u})-F^\infty(x,\textbf{u})}{\sum_{i=1}^m|u_i|^2+\sum_{i=1}^m|u_i|^{l_{2,i}+2}}=0.
$$
($H_5$): For all $x\in\mathbb{R}^N,\textbf{u}\in\mathbb{R}^m$, there exist $A',B'>0$ and $0<\beta'_i<l_{3,i}<\frac{4\alpha}{N}$ such that for all $1\leq i\leq m$, the function $F^\infty$ satisfies
\begin{align*}
&0\leq F^\infty(x,\textbf{u})\leq A'(\sum_{i=1}^m|u_i|^{\beta'_i+2}+\sum_{i=1}^m|u_i|^{l_{3,i}+2}),\\
&\partial_j F^\infty(x,\textbf{u})\leq B'(\sum_{i=1}^m|u_i|^{\beta'_i+1}+\sum_{i=1}^m|u_i|^{l_{3,i}+1}),\ \forall 1\leq j\leq m.
\end{align*}
($H_6$): For all $j \in\{0,1,\cdots,m\}$, we can write $F^\infty(x,u_1,\cdots,u_m) = F^\infty_{1j}(x,u_1,\cdots,u_m) + F^\infty_{2j}(x,u_1,\cdots,u_m) + F^\infty_{3j}(x,u_1,\cdots,u_{j-1},u_{j+1},\cdots,u_m)$. Furthermore, there exists $\sigma\in(0,+\infty)$ independent on $j$ such that for all $x\in \mathbb{R}^N,\textbf{u}\in\mathbb{R}^m$ and $\theta_j\geq 1$, we have:
	$$
	F^\infty_{1j}(x,u_1,u_2,\cdots,\theta_ju_j,\cdots,u_m)\geq\theta_j^{2+\sigma}F^\infty_{1j}(x,u_1,\cdots,u_m),
	$$
	$$
	F^\infty_{2j}(x,u_1,u_2,\cdots,\theta_ju_j,\cdots,u_m)\geq\theta_j^2F^\infty_{2j}(x,u_1,\cdots,u_m),
	$$
	there exists $L\in\mathbb{R}^N$ such that $$F^\infty_{3j}(x,u_1,\cdots,u_{j-1},u_{j+1},\cdots,u_m)=F^\infty_{3,j}(x+L,u_1,\cdots,u_{j-1},u_{j+1},\cdots,u_m)$$ for all $x\in\mathbb{R}^N$ and $\textbf{u}\in\mathbb{R}^m$, and there exists some $C > 0$ such that
	$$
	F^\infty(x,\textbf{u}) \leq C\sum_{j = 1}^mF^\infty_{1j}(x,\textbf{u})
	$$
	for all $x\in\mathbb{R}^N$ and $\textbf{u}\in\mathbb{R}^m$.
%($H_6$): There exists $\sigma\in(0,+\infty)$ such that for all $x\in\mathbb{R}^N, \textbf{u}\in\mathbb{R}^m$ and $\theta_i\geq 1,i=1,\ldots,m$, there holds
%$$
%F^\infty(x,\theta\textbf{u})\geq \theta_{max}^{\sigma+2}F^\infty(x,\textbf{u}),
%$$
%where $\theta_{max}=\max_{1\leq i\leq m}\theta_i$.

($H_7$): For all $x\in \mathbb{R}^N$ and $\textbf{u}\in \mathbb{R}^m,$ we have $F^\infty(x,\textbf{u})\leq F(x,\textbf{u})$ with strict inequality in a measurable set having a positive Lebesgue measure.\\

The class of nonlinearities satisfying $(H_1)-(H_7)$ includes many interesting functions. Here, we give some examples:
	\begin{itemize}
		\item[(1)] $m=2, N\geq1, \alpha > \frac{N}{2}, F(x,u_1,u_2)=p(x)|u_1|^2|u_2|^2+\mu_1|u_1|^{l_1}+\mu_2|u_2|^{l_2}+q(x)|u_1|^{k_1}|u_2|^{k_2}$ where $k_1,k_2 > 2, k_1+k_2<2+\frac{4\alpha}{N}, 2 < l_1, l_2 < 2+\frac{4\alpha}{N}, \mu_1, \mu_2 \geq 0$ and $p(x)=\frac{1}{1+|x|},q(x)=e^{-|x|}+1$,
		$$
		F^\infty(x,u_1,u_2)=\mu_1|u_1|^{l_1}+\mu_2|u_2|^{l_2}+|u_1|^{k_1}|u_2|^{k_2}.
		$$
		In this case,
		$$
		F_{11} = p(x)|u_1|^2|u_2|^2+\mu_1|u_1|^{l_1}+q(x)|u_1|^{k_1}|u_2|^{k_2}, F_{21} = \mu_2|u_2|^{l_2},
		$$
		$$
		F_{12} = p(x)|u_1|^2|u_2|^2+\mu_2|u_2|^{l_2}+q(x)|u_1|^{k_1}|u_2|^{k_2}, F_{22} = \mu_1|u_1|^{l_1},
		$$
		$$
		F_{11}^\infty = \mu_1|u_1|^{l_1} + |u_1|^{k_1}|u_2|^{k_2}, F_{21}^\infty = 0, F_{31}^\infty = \mu_2|u_2|^{l_2},
		$$
		$$
		F_{12}^\infty = \mu_2|u_2|^{l_2} + |u_1|^{k_1}|u_2|^{k_2}, F_{22}^\infty = 0, F_{32}^\infty = \mu_1|u_1|^{l_1}.
		$$
		
		\item[(2)] $m=2, N\geq1, \alpha > \frac{N}{2}, F(x,u_1,u_2)=\mu_1|u_1|^{4}+\mu_2|u_2|^{4}+q(x)|u_1|^{2}|u_2|^{2}$ where $q(x) =e^{-|x|}+1$, $\mu_1, \mu_2 > 0$ and
		$$
		F^\infty(x,u_1,u_2)=\mu_1|u_1|^{4}+\mu_2|u_2|^{4}+|u_1|^{2}|u_2|^{2}.
		$$
		In this case,
		$$
		F_{11} = \mu_1|u_1|^{4}+q(x)|u_1|^{2}|u_2|^{2}, F_{21} = \mu_2|u_2|^{4},
		$$
		$$
		F_{12} = \mu_2|u_2|^{4}+q(x)|u_1|^{2}|u_2|^{2}, F_{22} = \mu_1|u_1|^{4},
		$$
		$$
		F_{11}^\infty = \mu_1|u_1|^{4}, F_{21}^\infty = |u_1|^{2}|u_2|^{2}, F_{31}^\infty = \mu_2|u_2|^{4},
		$$
		$$
		F_{12}^\infty = \mu_2|u_2|^{4}, F_{22}^\infty = |u_1|^{2}|u_2|^{2}, F_{32}^\infty = \mu_1|u_1|^{4}.
		$$
		\item[(3)] $m=2, N\geq1, \alpha > \frac{3N}{4}, F(x,u_1,u_2)=\mu_1|u_1|^{l_1}+\mu_2|u_2|^{l_2}+q(x)\frac{|u_1|^{k_1}}{1+|u_1|}|u_2|^{k_2}$ where $k_1 >3, k_2 > 2, k_1+k_2<2+\frac{4\alpha}{N}, 2 < l_1, l_2 < 2+\frac{4\alpha}{N}, \mu_1, \mu_2 \geq 0$ and $q(x)=e^{-|x|}+1$,
		$$
		F^\infty(x,u_1,u_2)=\mu_1|u_1|^{l_1}+\mu_2|u_2|^{l_2}+\frac{|u_1|^{k_1}}{1+|u_1|}|u_2|^{k_2}.
		$$
		In this case,
		$$
		F_{11} = \mu_1|u_1|^{l_1}+q(x)\frac{|u_1|^{k_1}}{1+|u_1|}|u_2|^{k_2}, F_{21} = \mu_2|u_2|^{l_2},
		$$
		$$
		F_{12} = \mu_2|u_2|^{l_2}+q(x)\frac{|u_1|^{k_1}}{1+|u_1|}|u_2|^{k_2}, F_{22} = \mu_1|u_1|^{l_1},
		$$
		$$
		F_{11}^\infty = \mu_1|u_1|^{l_1} + \frac{|u_1|^{k_1}}{1+|u_1|}|u_2|^{k_2}, F_{21}^\infty = 0, F_{31}^\infty = \mu_2|u_2|^{l_2},
		$$
		$$
		F_{12}^\infty = \mu_2|u_2|^{l_2} + \frac{|u_1|^{k_1}}{1+|u_1|}|u_2|^{k_2}, F_{22}^\infty = 0, F_{32}^\infty = \mu_1|u_1|^{l_1}.
		$$
		\item[(4)]  $m=3, N\geq1, \alpha > N$,
		\begin{align*}
			F(x,u_1,u_2,u_3) =\sum_{j=1}^3\mu_j|u_j|^{l_j}+&\mu_{12}|u_1|^{k_1}|u_2|^{k_2}+\mu_{23}|u_2|^{k_3}|u_3|^{k_4}\\
			+&\mu_{13}|u_1|^{k_5}|u_3|^{k_6}+q(x)|u_1|^{t_1}|u_2|^{t_2}|u_3|^{t_3}
		\end{align*}
		where $k_i > 2$ for any $1 \leq i \leq 6$, $k_i+k_{i+1}<2+\frac{4\alpha}{N}$ for any $i = 1,3,5$, $2 < l_1, l_2, l_3 < 2+\frac{4\alpha}{N}, t_1, t_2, t_3 > 2, t_1+t_2+t_3 < 2+\frac{4\alpha}{N}, \mu_1, \mu_2, \mu_3, \mu_{12}, \mu_{13}, \mu_{23} \geq 0$ and $q(x) =e^{-|x|}+1$,
		\begin{align*}
		F^\infty(x,u_1,u_2,u_3) =\sum_{j=1}^3\mu_j|u_j|^{l_j}+&\mu_{12}|u_1|^{k_1}|u_2|^{k_2}+\mu_{23}|u_2|^{k_3}|u_3|^{k_4}\\
		+&\mu_{13}|u_1|^{k_5}|u_3|^{k_6}+|u_1|^{t_1}|u_2|^{t_2}|u_3|^{t_3}.
		\end{align*}
		In this case,
		$$
		F_{11} = \mu_1|u_1|^{l_1}+\mu_{12}|u_1|^{k_1}|u_2|^{k_2}+\mu_{13}|u_1|^{k_5}|u_3|^{k_6}+q(x)|u_1|^{t_1}|u_2|^{t_2}|u_3|^{t_3},
		$$
		$$
		F_{21}= \sum_{j=2}^3\mu_j|u_j|^{l_j}+\mu_{23}|u_2|^{k_3}|u_3|^{k_4},
		$$
		$$
		F_{12} = \mu_2|u_2|^{l_2}+\mu_{12}|u_1|^{k_1}|u_2|^{k_2}+\mu_{23}|u_2|^{k_3}|u_3|^{k_4}+q(x)|u_1|^{t_1}|u_2|^{t_2}|u_3|^{t_3},
		$$
		$$
		F_{22}= \mu_1|u_1|^{l_1}+\mu_3|u_3|^{l_3}+\mu_{13}|u_1|^{k_5}|u_3|^{k_6},
		$$
		$$
		F_{13} = \mu_3|u_3|^{l_3}+\mu_{23}|u_2|^{k_3}|u_3|^{k_4}+\mu_{13}|u_1|^{k_5}|u_3|^{k_6}+q(x)|u_1|^{t_1}|u_2|^{t_2}|u_3|^{t_3},
		$$
		$$
		F_{23}= \sum_{j=1}^2\mu_j|u_j|^{l_j}+\mu_{23}|u_2|^{k_3}|u_3|^{k_4},
		$$
		$$
		F^\infty_{11} = \mu_1|u_1|^{l_1}+\mu_{12}|u_1|^{k_1}|u_2|^{k_2}+\mu_{13}|u_1|^{k_5}|u_3|^{k_6}+|u_1|^{t_1}|u_2|^{t_2}|u_3|^{t_3},
		$$
		$$
		F^\infty_{21}=0, F^\infty_{31}= \Sigma_{j=2}^3\mu_j|u_j|^{l_j}+\mu_{23}|u_2|^{k_3}|u_3|^{k_4},
		$$
		$$
		F^\infty_{12} = \mu_2|u_2|^{l_2}+\mu_{12}|u_1|^{k_1}|u_2|^{k_2}+\mu_{23}|u_2|^{k_3}|u_3|^{k_4}+|u_1|^{t_1}|u_2|^{t_2}|u_3|^{t_3},
		$$
		$$
		F^\infty_{22}= 0, F^\infty_{32}= \mu_1|u_1|^{l_1}+\mu_3|u_3|^{l_3}+\mu_{13}|u_1|^{k_5}|u_3|^{k_6},
		$$
		$$
		F^\infty_{13} = \mu_3|u_3|^{l_3}+\mu_{23}|u_2|^{k_3}|u_3|^{k_4}+\mu_{13}|u_1|^{k_5}|u_3|^{k_6}+|u_1|^{t_1}|u_2|^{t_2}|u_3|^{t_3},
		$$
		$$
		F^\infty_{23}= 0, F^\infty_{33}= \Sigma_{j=1}^2\mu_j|u_j|^{l_j}+\mu_{23}|u_2|^{k_3}|u_3|^{k_4}.
		$$
	\end{itemize}
	
	All these particular examples satisfying all our hypotheses were not included in previous results \cite{LNW,Santosh 3,Santosh 4,GJ} that only addressed the case $\alpha = 1$, \cite{Santosh 1} for Hartree type equations, while\cite{Santosh 2} studied $0<\alpha<1$ of Choquard type. In these papers, specific nonlinearities were studied and most cases are specific ones of our general nonlinearities. For $F(u_1,u_2)$ such that $\partial_iF(u_1,u_2)$ is nondecreasing for $i = 1,2$ and $F$ satisfies other suitable conditions, the readers can refer to \cite{S}.
%\begin{itemize}
%	\item[\cite{LNW}] $m\geq 2, \alpha = 1, N = 1, f(x,\textbf{u}) = \frac1p\Sigma_{i,j=1}^mb_{ij}|u_iu_j|^p$ where $2 \leq p < 3$ and $b_{ij} > 0$;
%	\item[\cite{Santosh 3}] $m = 2, \alpha = 1, N = 1, f(x,\textbf{u}) = a|u_1|^p + b|u_2|^r + c|u_1u_2|^q$ where $2 < p,r,2q < 6$ and $a, b, c > 0$;
%	\item[\cite{Santosh 4}] $m = 3, \alpha = 1, N = 1, f(x,\textbf{u}) = \frac1p\Sigma_{i,j=1}^3b_{ij}|u_iu_j|^p$ where $2 \leq p < 3$ and $b_{ij} > 0$;
%	\item[\cite{GJ}] $m = 2, \alpha = 1, N \geq 1, f(x,\textbf{u}) = a|u_1|^{p_1} + b|u_2|^{p_2} + c|u_1|^{r_1}|u_2|^{r_2}$ where $2 < p_1, p_2 < 2 + \frac4N, r_1, r_2 > 1, r_1+r_2 < 2 + \frac4N$ and $a, b, c > 0$;
%	\item[\cite{Santosh 1}] $m =2, 3, \alpha =1, N \geq 1,$
%	\begin{align*}
%	J(\textbf{u})=\frac{1}{2}\int_{\mathbb{R}^N}\sum_{i=1}^m|\Lambda^\alpha u_i|^2dx -\frac{1}{2p}\Sigma_{i,j=1}^m\int_{\mathbb{R}^N}(W \star |u_i|^p)|u_j|^pdx
%	\end{align*}
%	where $2 \leq p < \frac{2r-1}{r} + \frac{2}{N}$ with $\frac{1}{r} < \frac{2}{N}$ and the convolution potential $W$ satisfies suitable conditions;
%	\item[\cite{Santosh 2}] $m = 2, 0 < \alpha <1, N \geq 2,$
%	\begin{align*}
%		& J(\textbf{u})=\frac{1}{2}\int_{\mathbb{R}^N}\sum_{i=1}^2|\Lambda^\alpha u_i|^2dx\\
%		-&\int\int_{\mathbb{R}^N}\left(\frac{a}{2p_1}\frac{|u_1(x)|^{p_1}|u_1(y)|^{p_1}}{|x-y|^{N-\beta}} + \frac{b}{2p_2}\frac{|u_2(x)|^{p_2}|u_2(y)|^{p_2}}{|x-y|^{N-\beta}} + \frac{c}{q}\frac{|u_1(x)|^{q}|u_2(y)|^{q}}{|x-y|^{N-\beta}}\right)dxdy
%	\end{align*}
%	where $0 < \beta < N$, $2 \leq p_1, p_2, q < \frac{N+2\alpha+\beta}{N}$ and $a, b, c > 0$.
%\end{itemize}
For more related topics, the readers can see \cite{GLN} for the constraint $\left\langle u_k,u_j\right\rangle = \delta_{k,j}$, \cite{BZZ} for the $L^2$ supercritical case using a new approach based on bifurcation theory and the continuation method, \cite{BS} for infinitely many solutions in the $L^2$ supercritical case, \cite{BCMP} for the energy critical case in a fractional setting, \cite{LZ} for a mixed situation ($L^2$ subcritical and $L^2$ supercritical), and references therein.

%$$
%F(x,u_1,u_2)=p(x)(|u_1|^2+|u_2|^2)+q(x)|u_1|^{k_1}|u_2|^{k_2},\quad k_1,k_2<\frac{4\alpha}{3}+1
%$$
%where $p(x)=\frac{1}{1+|x|},q(x)=e^{-|x|}+1$.
%$$
%F^\infty(x,u_1,u_2)=|u_1|^{k_1}|u_2|^{k_2}
%$$

Now, we state the main result of this paper.
\begin{theorem}\label{th1.1}
Suppose that $(H_1)-(H_7)$ hold. Then there exists $\textbf{u}\in S_c$ such that $J(\textbf{u})=I_{c_1,\ldots,c_m}.$ $\textbf{u}$ is called a normalized ground state solution (NGSS) of the system \eqref{1.1}.
\end{theorem}

From our arguments, an intermediate result is obtained.
\begin{theorem}\label{th1.2}
Suppose that $(H_2), (H_5), (H_6)$ hold true for $F^\infty$. Then there exists $\textbf{u}\in S_c$ such that $J^\infty(\textbf{u})=I^\infty_{c_1,\ldots,c_m}$ where
\begin{align*}
&J^\infty(\textbf{u})=\frac{1}{2}\int\sum_{i=1}^m|\Lambda^\alpha u_i|^2-\int F^\infty(x,\textbf{u})dx,\\
&I^\infty_{c_1,\ldots,c_m}=\inf\{J^\infty(\textbf{u}),\textbf{u}\in S_c\}.
\end{align*}
\end{theorem}

\begin{remark}
	(1) To the best of our knowledge, this is the first work concerning \textbf{Fractional Schr\"{o}dinger System with General Nonlinearities} allowing $\alpha > 1$.
	
	(2) Our discussions are for local nonlinearities, but can be easily extended to Hartree type nonlinearities (and combination of both). Particularly, situations in \cite{Santosh 1,Santosh 2} can be addressed.
	
	(3) We don't used rearrangement techniques like \cite{S,GJ}. Hence it not necessary to assume that $F(x,{\bf u})$ and $F^\infty(x,{\bf u})$ are radial w.r.t. $x \in \mathbb{R}^N$.

(4) By establishing the uniqueness of the solution to the Cauchy problem associated to \eqref{1.1}, we can easily show the orbital stability of standing waves by following the steps in \cite{HS}.

(5) The situation with orthonormal functions $u_i$ in the constrained is much more challenging than the constrained $\int u_1^2+u_2^2+\cdots+u_n^2=c,$ where things are reduced to the one equation setting, see \cite{H}.
\end{remark}

The following improved concentration-compactness principle plays a crucial role in our proofs.
\begin{lemma}\label{le1.3}
 Assume that the vectorial sequence $\textbf{u}_{n}=(u_{n,1},\ldots,u_{n,m})\subset H^{\alpha}$ satisfies
$$
||u_{n,i}||^2_{L^{2}}=\mu_i>0,\ i=1,\ldots,m,
$$
$$
\sum_{i=1}^m||u_{n,i}||_{H^{\alpha}}<M,\quad \forall n\in\mathbb{N^+}.
$$
Then, there exists a vectorial subsequence $\textbf{u}_{n}$ (that we still denote by $\textbf{u}_{n}$) for which one of the following properties
holds.
\newline
(1)Compactness: there exists a sequence $\{y_{n}\}_{n\in\mathbb{N}^+}\subset\mathbb{R}^N$ such that for any $\epsilon>0$, there exists $0<r<\infty$ with
$$
\sum_{i=1}^m\int_{|x-y_{n}|\leq r}|u_{n,i}(x)|^{2}dx\geq \mu-\epsilon,
$$
where $\mu=\sum_{i=1}^{m}\mu_i.$
\newline
(2)Vanishing: For all $r<\infty$, it follows that
$$
\lim_{n\rightarrow\infty}\sup_{y\in \mathbb{R}^N}\sum_{i=1}^m\int_{|x-y|\leq r}|u_{n,i}|^{2}dx=0.
$$
(3)Dichotomy: There exists  constants $\beta_i\in(0,\mu_i), i=1,\ldots,m$ and two bounded vectorial sequences $\textbf{v}_{n}=(v_{n,1},\ldots,v_{n,m}),\textbf{w}_{n}=(w_{n,1},\ldots,w_{n,m})\subset H^{\alpha}$ such that for $1\leq i,j\leq m$, we have
$$
\mbox{supp}\ v_{n,i}\bigcap \mbox{supp}\ w_{n,j}=\emptyset,
$$
$$
|v_{n,i}|+|w_{n,i}|\leq |u_{n,i}|.
$$
$$
||v_{n,i}||^{2}_{L^{2}}\rightarrow \beta_i,\quad ||w_{n,i}||^{2}_{L^{2}}\rightarrow (\mu_i-\beta_i),\  as\ n\rightarrow\infty.
$$
$$
||u_{n,i}-v_{n,i}-w_{n,i}||_{L^{p}}\rightarrow 0,\ \textrm{for}\ 2\leq p<2^*_\alpha,
$$
where 
\begin{align*}
2^*_\alpha=\left\{
      \begin{array}{ll}
        \infty, & N \leq 2\alpha, \\
        \frac{2N}{N-2\alpha}, & N > 2\alpha.
      \end{array}
    \right.
\end{align*} 
\begin{eqnarray*}
\liminf_{n\rightarrow\infty}\{\langle(-\triangle)^{\alpha}u_{n,i},u_{n,i}\rangle-\langle(-\triangle)^{\alpha}v_{n,i},v_{n,i}\rangle-\langle(-\triangle)^{\alpha}w_{n,i},w_{n,i}\rangle\}\geq 0.
\end{eqnarray*}
\end{lemma}

This result is obtained by considering the L$\mathrm{\acute{e}}$vy concentration functions:
\begin{align*}
Q_n(r)=\sup_{y\in\mathbb{R}^N}\int_{|x-y|<r}\sum_{i=1}^m|u_{n,i}|^2dx.
\end{align*}
%Following the approach of Proposition 1.7.6 \cite{Ca}, we can derive the compactness and vanishing property.% The dichotomy is obtained with a similar argument of Lemma 2.4 \cite{Feng}. Here, we omit the details for short.

\section{A few technical Lemmata}

To investigate the constrained variational problems, we firstly give some estimates.
\begin{lemma}\label{le2.1}
Suppose that $(H_1)$ and $(H_5)$ holds true for $F$ and $F^\infty$ respectively. Then, we have the following conclusions:
\newline
(1) $J(\textbf{u})$ and $J^\infty(\textbf{u})$ are  $C^1$ functional satisfying
\begin{align*}
&||J(\textbf{u})||_{H^{-1}}\leq C\sum_{i=1}^m(||u_i||_{H^\alpha}+||u_i||^{1+\frac{4\alpha}{N}}_{H^\alpha}),\\
&||J^\infty(\textbf{u})||_{H^{-1}}\leq C\sum_{i=1}^m(||u_i||_{H^\alpha}+||u_i||^{1+\frac{4\alpha}{N}}_{H^\alpha}),\ \forall\ \textbf{u}\in H^\alpha.
\end{align*}
\newline
(2) There exist two positive constants $C_{A,m,c_1,\ldots,c_m}$ and $C_{A',m,c_1,\ldots,c_m}$ such that
\begin{align*}
&J(\textbf{u})\geq \frac{1}{4}\sum_{i=1}^m||u_i||_{\dot{H}^\alpha}^2-C_{A,m,c_1,\ldots,c_m},\\
&J^\infty(\textbf{u})\geq \frac{1}{4}\sum_{i=1}^m||u_i||_{\dot{H}^\alpha}^2-C_{A',m,c_1,\ldots,c_m},\ \forall\ \textbf{u}\in H^\alpha.
\end{align*}
\newline
(3) $I_{c_1,\ldots,c_m}>-\infty$ and $I^\infty_{c_1,\ldots,c_m}>-\infty$. In addition, any minimizing sequence of $I_{c_1,\ldots,c_m}$ and $I^\infty_{c_1,\ldots,c_m}$ is bounded in $H^\alpha$.
\newline
(4) The mappings $(c_1,\ldots,c_m)\mapsto I_{c_1,\ldots,c_m}$ and $(c_1,\ldots,c_m)\mapsto I^\infty_{c_1,\ldots,c_m}$ are continuous on $(0,+\infty)^m$.

\end{lemma}
\paragraph{Proof}

We only prove the assertion for $J(\textbf{u})$ and $I_{c_1,\ldots,c_m}$ here. The proofs for $J^\infty(\textbf{u})$ and $I^\infty_{c_1,\ldots,c_m}$ are similar.

First, we introduce the cut off function $\phi(x)\in C^\infty_0( \mathbb{R},\mathbb{R}^+)$ such that
\begin{align*}
\phi(x)=\left\{
                   \begin{array}{ll}
                     1, & |x|\leq 1, \\
                     0,&|x|\geq 2.
                   \end{array}
                 \right.
\end{align*}
By direct calculations, there exists a constant $K>0$ such that for $1\leq i\leq m$, we have
\begin{align*}
&||\prod_{j=1}^m\phi(u_j)\partial_i F(x,\textbf{u})||_2\leq K\sum_{j=1}^m||u_j||_2,\\
&||(1-\prod_{j=1}^m\phi(u_j))\partial_i F(x,\textbf{u})||_p\leq K\sum_{j=1}^m||u_j||_q^{1+\frac{4\alpha}{N}},
\end{align*}
where
\begin{align*}
p=\left\{
    \begin{array}{ll}
      \frac{2N}{N+2\alpha}, & N\geq 2, \\
      \frac{4}{3}, & N=1,
    \end{array}
  \right.
\ q=\big(\frac{4\alpha}{N}+1\big)p.
\end{align*}

For any $\psi\in H^\alpha$, we have
\begin{align*}
\langle \partial_i F(x,\textbf{u}),\psi\rangle&=\langle \prod_{j=1}^m\phi(u_j) \partial_i F(x,\textbf{u}),\psi\rangle+\langle (1-\prod_{j=1}^m\phi(u_j))\partial_i F(x,\textbf{u}),\psi\rangle\\
&\leq ||\prod_{j=1}^m\phi(u_j) \partial_i F(x,\textbf{u})||_2||\psi||_2+||(1-\prod_{j=1}^m\phi(u_j))\partial_i F(x,\textbf{u})||_p||\psi||_{p^*}\\
&\leq K\sum_{i=1}^m\Big(||u_i||_2||\psi||_2+||u_i||_q^{1+\frac{4\alpha}{N}}||\psi||_{p^*}\Big)\\
&\leq C_K\sum_{i=1}^m(||u_i||_{H^\alpha}+||u_i||^{1+\frac{4\alpha}{N}}_{H^\alpha})||\psi||_{H^\alpha},
\end{align*}
where $p^*=\frac{p}{p-1}$. The last inequality is obtained by using Sobolev embedding theorem.
Thus, we get
\begin{align}\label{2.1}
||\partial_i F(x,\textbf{u})||_{H^{-\alpha}}\leq C_K\sum_{i=1}^m(||u_i||_{H^\alpha}+||u_i||^{1+\frac{4\alpha}{N}}_{H^\alpha}).
\end{align}

Note that
\begin{align*}
\int F(x,\textbf{u})dx\leq A\sum_{i=1}^m(||u_i||_2^2+||u_i||^{l_{1,i}+2}_{l_{1,i}+2})\leq C_A\sum_{i=1}^m(||u_i||^2_{H^\alpha}+||u_i||^{l_{1,i}+2}_{H^\alpha}).
\end{align*}
It implies that $J(\textbf{u})$ is $C^1$ on $H^\alpha$ and
\begin{align*}
J'(\textbf{u})\textbf{v}=\int \sum_{i=1}^{m}\Lambda^\alpha u_i\Lambda^\alpha v_idx-\int \sum_{i=1}^m\partial_iF(x,\textbf{u})v_idx, \ \forall \ \textbf{v}:=(v_1,\ldots,v_m)\in H^\alpha.
\end{align*}
From \eqref{2.1}, we have
\begin{align}\label{2.2}
||J'(\textbf{u})||_{H^{-\alpha}}\leq  C_K\sum_{i=1}^m(||u_i||_{H^\alpha}+||u_i||^{1+\frac{4\alpha}{N}}_{H^\alpha}),\ \forall \ \textbf{u}\in H^\alpha.
\end{align}

Next, we estimate $J(\textbf{u})$.
 For any $\textbf{u}\in S_c$, there holds
\begin{align*}
\int F(x,\textbf{u})dx\leq Ac+A\sum_{i=1}^m||u_i||_{l_{1,i}+2}^{l_{1,i}+2}.
\end{align*}
From the fractional Gagliardo-Nirenberg inequality \cite{HMWO} and Young's inequality, for $1\leq i\leq m$, we have
\begin{align*}
||u_i||_{l_{1,i}+2}^{l_{1,i}+2}&\leq C||u_i||^{\frac{2\alpha(l_{1,i}+2)-Nl_{1,i}}{2\alpha}}_2||u_i||^{\frac{Nl_{1,i}}{2\alpha}}_{\dot{H}^\alpha}\\
&\leq \frac{1}{4Am}(||u_i||_{\dot{H}^\alpha}^{\frac{Nl_{1,i}}{2\alpha}})^{\frac{4\alpha}{Nl_{1,i}}}+C(\frac{1}{4AmC})^{\frac{Nl_{1,i}}{Nl_{1,i}-4\alpha}}
(||u_i||_2^{\frac{2\alpha(l_{1,i}+2)-Nl_{1,i}}{2\alpha}})^{\frac{4\alpha}{4\alpha-Nl_{1,i}}}\\
&\leq \frac{1}{4Am}||u_i||_{\dot{H}^\alpha}^2+C_{A,m}||u_i||_2^{\frac{4\alpha(l_{1,i}+2)-2Nl_{1,i}}{4\alpha-Nl_{1,i}}}.
\end{align*}
Therefore,
\begin{align*}
J(\textbf{u})\geq \frac{1}{4}\sum_{i=1}^m||u_i||_{\dot{H}^\alpha}^2-Ac-C_{A,m,c_1,\ldots,c_m}.
\end{align*}
It implies that $J(\textbf{u})$ is bounded from below on $S_c$ and any minimizing sequence of the constrained variational problem $I_{c_1,\ldots,c_m}$ is bounded in $H^\alpha$.

Now, we prove the continuity of the map $(c_1,\ldots,c_m)\mapsto I_{c_1,\ldots,c_m}$ on $(0,\infty)^m$. For any fixed $c=(c_1,\ldots,c_m)$, we consider a sequence $c_n=(c_{n,1},\ldots,c_{n,m})$ satisfying $c_{n,i}\rightarrow c_i$ for $1\leq i\leq m$. Then, for any $n\in\mathbb{N}^+$, there exists $\textbf{u}_n=(u_{n,1},\ldots,u_{n,m})\in S_{c_n}$ with $\int |u_{n,i}|^2dx=c_{n,i},i=1,\ldots,m$ such that
\begin{align*}
I_{c_{n,1},\ldots,c_{n,m}}\leq J(\textbf{u}_n)\leq I_{c_{n,1},\ldots,c_{n,m}}+\frac{1}{n}.
\end{align*}

Using a similar argument as the previous one, there exists a constant $K>0$ such that
\begin{align*}
\sum_{i=1}^m||u_{n,i}||_{H^\alpha}\leq K,\ \forall\ n\in\mathbb{N}^+.
\end{align*}
 Let $\tilde{\textbf{u}}_n=(\tilde{u}_{n,1},\ldots,\tilde{u}_{n,m})$ satisfy $\tilde{u}_{n,i}=\sqrt{\frac{c_i}{c_{n,i}}}u_{n,i},i=1,\ldots,m$. It is obvious that $\tilde{\textbf{u}}_n\in S_c$. Moreover, there exists $n_1>0$ such that
\begin{align*}
||\textbf{u}_n-\tilde{\textbf{u}}_n||_{H^\alpha}&:=\sum_{i=1}^m||u_{n,i}-\tilde{u}_{n,i}||_{H^\alpha}\\
&\leq \sum_{i=1}^m\big|\sqrt{\frac{c_i}{c_{n,i}}}-1\big|||u_{n,i}||_{H^\alpha}\\
&\leq K+1,\ \forall \ n\geq n_1.
\end{align*}

Due to the estimate \eqref{2.2}, we have:
\begin{align*}
|J(\textbf{u}_n)-J(\tilde{\textbf{u}}_n)|&\leq \sup_{||\textbf{u}||_{H^\alpha}\leq 2K+1}||J'(\textbf{u})||_{H^{-\alpha}}||\textbf{u}_n-\tilde{\textbf{u}}_n||_{H^\alpha}\\
&\leq C_K\sum_{i=1}^m|1-\sqrt{\frac{c_i}{c_{n,i}}}|.
\end{align*}
Therefore, we have
\begin{align*}
I_{c_{n,1},\ldots,c_{n,m}}\geq J(\textbf{u}_n)-\frac{1}{n}\geq J(\tilde{\textbf{u}}_n)+C_K\sum_{i=1}^m|1-\sqrt{\frac{c_i}{c_{n,i}}}|-\frac{1}{n}
\end{align*}
which implies that
\begin{align}\label{2.3}
\liminf_{n\rightarrow\infty}I_{c_{n,1},\ldots,c_{n,m}}\geq I_{c_1,\ldots,c_m}.
\end{align}

On the other hand, let $\textbf{u}_n\in S_c$ be a minimizing sequence satisfying $\lim_{n\rightarrow\infty}J(\textbf{u}_n)=I_{c_1,\ldots,c_m}$ and $\sum_{i=1}^m||u_{n,i}||_{H^\alpha}\leq K$. Then, for any  sequence $c_n=(c_{n,1},\ldots,c_{n,i})$ with $c_{n,i}\rightarrow c_i$, $i=1,\ldots,m$, we set $\tilde{u}_{n,i}=\sqrt{\frac{c_{n,i}}{c_i}}u_{n,i}$. Therefores
\begin{align*}
I_{c_{n,1},\ldots,c_{n,m}}\leq J(\tilde{\textbf{u}}_n)\leq J(\textbf{u}_n)+C_K\sum_{i=1}^m|1-\sqrt{\frac{c_{n,i}}{c_i}}|.
\end{align*}
Thus, $\limsup_{n\rightarrow\infty}I_{c_{n,1},\ldots,c_{n,m}}\leq I_{c_1,\ldots,c_m}$. Combining this and \eqref{2.3}, we get
\begin{align*}
\lim_{n\rightarrow\infty}I_{c_{n,1},\ldots,c_{n,m}}=I_{c_1,\ldots,c_m}
\end{align*}
which completes the proof.
$\square$

\begin{lemma}\label{le2.2}
$~$
\begin{enumerate}
  \item [(1)] Assume that $(H_1)$ and $(H_2)$ hold true for $F$. Then, $I_{c_1,\ldots,c_m}<0$.
  \item [(2)] Assume that $(H_2)$ and $(H_5)$ hold true for $F^\infty$. Then, $I^\infty_{c_1,\ldots,c_m}<0$.
\end{enumerate}
\end{lemma}
\paragraph{Proof}
 Let $\phi$ be a fixed radial function with $||\phi||_2=1$. Define
\begin{align*}
\Phi_\lambda=\lambda^{\frac{N}{2}}\Phi(\lambda x)=\lambda^{\frac{N}{2}}(\phi_1(\lambda x),\ldots,\phi_m(\lambda x)),
\end{align*}
where $\phi_i=\sqrt{c_i}\phi$. $0<\lambda<1$ is a small constant to be chosen later.

One can see that
\begin{align*}
J(\Phi_\lambda)&=\frac{1}{2}||\lambda^{\frac{N}{2}}\Phi_\lambda||^2_{\dot{H}^\alpha}-\int F(x,\lambda^{\frac{N}{2}}\phi_1(\lambda x),\ldots,\lambda^{\frac{N}{2}}\phi_m(\lambda x))dx\\
&\leq \frac{1}{2}\int \sum_{i=1}^m|\xi|^{2\alpha}|\widehat{\lambda^{\frac{N}{2}}\phi_i(\lambda x)}|^2d\xi-\int_{|x|\geq R} F(x,\lambda^{\frac{N}{2}}\phi_1(\lambda x),\ldots,\lambda^{\frac{N}{2}}\phi_m(\lambda x))dx\\
&\leq \frac{\lambda^{2\alpha}}{2}\int |\frac{\xi}{\lambda}|^{2\alpha}|\widehat{\phi}\big(\frac{\xi}{\lambda}\big)|^2d\big(\frac{\xi}{\lambda}\big)
-\lambda^{\frac{Ns}{2}}\beta\int_{|x|\geq R}|x|^{-t}\phi_1^{s_1}(\lambda x)\cdots\phi_m^{s_m}(\lambda x)dx\\
&\leq \frac{\lambda^{2\alpha}}{2}||\Phi||^2_2-\beta\lambda^{\frac{Ns}{2}-N+t}\int_{|y|\geq \lambda R}|y|^{-t}\phi_1^{s_1}(y)\cdots\phi_m^{s_m}(y)dy
\end{align*}
Thus, there exist two positive constants $C_1$ and $C_2$ such that
\begin{align*}
J(\Phi_\lambda)\leq \lambda^{2\alpha}(C_1-\lambda^{\frac{Ns}{2}-N+t-2\alpha}C_2).
\end{align*}
Note that $\frac{Ns}{2}-N+t-2\alpha<0$. Therefore, it completes the proof by choosing $\lambda<\Big(\frac{C_2}{C_1}\Big)^{-\frac{1}{\frac{Ns}{2}-N+t-2\alpha}}$.~~~~~~~~~ $\square$

Now, we prove the sub-additivity of $I_{c_1,\ldots,c_m}$ and $I^\infty_{c_1,\ldots,c_m}$.
\begin{lemma}\label{le2.3}(Strict sub-additivity)
$~$
\begin{enumerate}
  \item [(1)] Assume that $(H_1)$, $(H_2)$ and $(H_3)$ hold true for $F$. Then, for any $a_i\in (0,c_i), i=1,\ldots,m$, we have $I_{c_1,\ldots,c_m}\leq I_{a_1,\ldots,a_m}+I_{c_1-a_1,\ldots,c_m-a_m}$.
\newline
\item [(2)] Assume that $(H_2)$, $(H_5)$ and $(H_6)$ hold true for $F^\infty.$ Then, for any $a_i\in (0,c_i), i=1,\ldots,m$, we have $I^\infty_{c_1,\ldots,c_m}< I^\infty_{a_1,\ldots,a_m}+I^\infty_{c_1-a_1,\ldots,c_m-a_m}$.
    \end{enumerate}
\end{lemma}
\paragraph{Proof} First, we will prove the assertion (1). Let $\textbf{u}_n = (u_{n,1},\cdots,u_{n,m}) \in S_{a}$, $\textbf{v}_n = (v_{n,1},\cdots,v_{n,m}) \in S_{c-a}$ be bounded minimizing sequences of $I_{a_1,\ldots,a_m}$ and $I_{c_1-a_1,\ldots,c_m-a_m}$ respectively, satisfying $\lim_{n\rightarrow\infty}J(\textbf{u}_n)=I_{a_1,\ldots,a_m}, \lim_{n\rightarrow\infty}J(\textbf{v}_n)=I_{c_1-a_1,\ldots,c_m-a_m}$. By a density argument, we may suppose that $u_{n,i}, v_{n,i}$ have compact support for all $1 \leq i \leq m$. For $\textbf{u} = (u_1, \cdots, u_m) \in H^\alpha$, we define
$$
J_1(\textbf{u}) := \frac{1}{2}\int|\Lambda^\alpha u_1|^2-\int F_{11}(x,u_1,\cdots,u_{m})dx,
$$
$$
J_2(u_2,\cdots,u_m) := \sum_{i=2}^m\frac{1}{2}\int|\Lambda^\alpha u_i|^2-\int F_{21}(x,u_2,\cdots,u_m)dx,
$$
where $F_{11}$ and $F_{21}$ are given by $(H_3)$. Set
$$
M_1 := \frac{1}{a_1}\lim_{n \to \infty}J_1(\textbf{u}_n),
$$
$$
M_2 := \frac{1}{c_1-a_1}\lim_{n \to \infty}J_1(\textbf{v}_n).
$$
WLOG, we assume that $M_1 \leq M_2$. We denote $\tilde{v}_{n,i}(x) = v_{n,i}(x+b_iL)$ for all $2 \leq i \leq m$ where $L$ is given by $(H_3)$. For any $\epsilon > 0$, by chosing $b_i$ large enough, direct calculations yield that
\begin{align*}
&\lim_{n \to \infty}J_2(u_{n,2}+\tilde{v}_{n,2},\cdots,u_{n,m}+\tilde{v}_{n,m}) \\
&\leq \lim_{n \to \infty}J_2(u_{n,2},\cdots,u_{n,m}) + \lim_{n \to \infty}J_2(\tilde{v}_{n,2},\cdots,\tilde{v}_{n,m}) + \epsilon \\
&= \lim_{n \to \infty}J_2(u_{n,2},\cdots,u_{n,m}) + \lim_{n \to \infty}J_2(v_{n,2},\cdots,v_{n,m}) + \epsilon.
\end{align*}
Moreover, we can choose suitably $b_i$ and $R_n$ such that $\cup_{i \geq 1}\mbox{supp}\ u_{n,i} \subset B_{R_n}$ and $\cup_{i \geq 2}\mbox{supp}\ v_{n,i} \subset B_{R_n}^c$. Then by $(H_3)$ we have
\begin{align*}
&\lim_{n \to \infty}J_1(\sqrt{\frac{c_1}{a_1}}u_{n,1},u_{n,2}+\tilde{v}_{n,2},\cdots,u_{n,m}+\tilde{v}_{n,m}) \\
&\leq \lim_{n \to \infty}\frac{c_1}{a_1}J_1(u_{n,1},\cdots,u_{n,m}) \\
&\leq a_1M_1 + (c_1-a_1)M_2.
\end{align*}
Noticing that $(\sqrt{\frac{c_1}{a_1}}u_{n,1},u_{n,2}+\tilde{v}_{n,2},\cdots,u_{n,m}+\tilde{v}_{n,m}) \in S_c$ we deduce that
\begin{align*}
& I_{c_1,\cdots,c_m}\\
\leq&\lim_{n \to \infty}J(\sqrt{\frac{c_1}{a_1}}u_{n,1},u_{n,2}+\tilde{v}_{n,2},\cdots,u_{n,m}+\tilde{v}_{n,m}) \\
=& \lim_{n \to \infty}(J_1(\sqrt{\frac{c_1}{a_1}}u_{n,1},u_{n,2}+\tilde{v}_{n,2},\cdots,u_{n,m}+\tilde{v}_{n,m}) + J_2(u_{n,2}+\tilde{v}_{n,2},\cdots,u_{n,m}+\tilde{v}_{n,m})) \\
\leq& \lim_{n \to \infty}J_1(\textbf{u}_n) + \lim_{n \to \infty}J_1(\textbf{v}_n) + \lim_{n \to \infty}J_2(u_{n,2},\cdots,u_{n,m}) + \lim_{n \to \infty}J_2(v_{n,2},\cdots,v_{n,m}) + \epsilon \\
=& \lim_{n \to \infty}J(\textbf{u}_n) + \lim_{n \to \infty}J(\textbf{v}_n) + \epsilon.
\end{align*}
By the arbitrariness of $\epsilon$ we complete the proof of (i).

Now, we focus on studying $I^\infty_{c_1,\ldots,c_m}$. Let $\textbf{u}_n \in S_a, \textbf{v}_n\in S_{c-a}$ be two bounded minimizing sequences of $I^\infty_{a_1,\ldots,a_m}$ and $I^\infty_{c_1-a_1,\ldots,c_m-a_m}$ respectively. By $(H_6)$, similar to the proof of (i), we know that $I^\infty_{c_1,\ldots,c_m} \leq I^\infty_{a_1,\ldots,a_m}+I^\infty_{c_1-a_1,\ldots,c_m-a_m}$. Furthermore, the equality infers that
$$
 \lim_{n \to \infty}F^\infty_{1j}(x,\textbf{u}_n) = 0, \forall j \in \{1,\cdots,m\},
$$
implying that $\lim_{n \to \infty}F^\infty(x,\textbf{u}_n) = 0$. This yields that $I^\infty_{a_1,\cdots,a_m} \geq 0$, in a contradiction with Lemma \ref{le2.2} (2).
$\square$

\begin{lemma}\label{le2.4}
 Assume that $F$ and $F^\infty$ satisfy $(H_1)$, $(H_2)$, $(H_3)$  and $(H_2)$, $(H_5)$, $(H_6)$, $(H_7)$ respectively. Then, it follows that $I_{c_1,\ldots,c_m}\leq I^\infty_{c_1,\ldots,c_m}$ for any $c_i, i=1,\ldots,m$. Moreover, for any $a_i\in (0,c_i), i=1,\ldots,m$, we have $I_{c_1,\ldots,c_m}< I_{a_1,\ldots,a_m}+I^\infty_{c_1-a_1,\ldots,c_m-a_m}$.
\end{lemma}
\paragraph{Proof}
The result is obtained by Lemma \ref{le2.3} and hypothesis $(H_7)$ immediately.
$\square$

\section{Proof of Theorem \ref{1.2}}
In this section, we will prove Theorem \ref{1.2} by the improved concentration-compactness principle (Lemma \ref{le1.3}). Assume that $\textbf{u}_n$ is a minimizing sequence of the constrained variational problems $I^\infty_{c_1,\ldots,c_m}$. First, we will rule out vanishing by contradiction. For that, we suppose that $\textbf{u}_n$  is vanishing. That is,  for all $r<\infty$, it follows that
$$
\lim_{n\rightarrow\infty}\sup_{y\in \mathbb{R}^N}\sum_{i=1}^m\int_{|x-y|\leq r}|u_{n,i}|^{2}dx=0.
$$
From Lemma 2.4 in \cite{Guo}, one can see $\lim_{n\rightarrow\infty}\sum_{i=1}^m||u_{n,i}||_p=0$ for $p\in(2,2^*_\alpha)$.

Recalling the assumption $(H_5)$, we have
\begin{align*}
\lim_{n\rightarrow\infty}\int F^\infty(x,\textbf{u})dx=0.
\end{align*}
It obtains that $\liminf_{n\rightarrow\infty}J^\infty(\textbf{u}_n)\geq 0$ which contradicts with Lemma \ref{le2.2}.

Next, we will prove that the dichotomy does not occur by contradiction. Suppose that there exists constants $a_i\in(0,c_i)$ and two bounded vectorial sequences $\textbf{v}_{n}=(v_{n,1},\ldots,v_{n,m}),\textbf{w}_{n}=(w_{n,1},\ldots,w_{n,m})\subset H^{\alpha}$ such that for $1\leq i,j\leq m$, we have
$$
\mbox{supp}\ v_{n,i}\bigcap \mbox{supp}\ w_{n,j}=\emptyset,
$$
$$
|v_{n,i}|+|w_{n,i}|\leq |u_{n,i}|.
$$
$$
||v_{n,i}||^{2}_{L^{2}}:=a_{n,i}\rightarrow a_i,\quad ||w_{n,i}||^{2}_{L^{2}}:=b_{n,i}\rightarrow (c_i-a_i),\  as\ n\rightarrow\infty.
$$
$$
||u_{n,i}-v_{n,i}-w_{n,i}||_{L^{p}}\rightarrow 0,\ \textrm{for}\ 2\leq p<2^*_\alpha,
$$
where 
\begin{align*}
2^*_\alpha=\left\{
      \begin{array}{ll}
        \infty, & N\leq 2\alpha \\
        \frac{2N}{N-2\alpha}, & N> 2\alpha.
      \end{array}
    \right.
\end{align*}
\begin{eqnarray*}
\liminf_{n\rightarrow\infty}\{\langle(-\triangle)^{\alpha}u_{n,i},u_{n,i}\rangle-\langle(-\triangle)^{\alpha}v_{n,i},v_{n,i}\rangle-\langle(-\triangle)^{\alpha}w_{n,i},w_{n,i}\rangle\}\geq 0.
\end{eqnarray*}

From the arguments in Lemma \ref{2.1}, there exists a constant $C_K>0$ such that
\begin{align*}
&|\int(F^\infty(x,\textbf{u}_n)-F^\infty(x,\textbf{v}_n)-F^\infty(x,\textbf{w}_n))dx|\\
\leq&|\int(F^\infty(x,\textbf{u}_n)-F^\infty(x,\textbf{v}_n+\textbf{w}_n)dx|\\
\leq&\sup_{||\textbf{u}||_{H^\alpha}\leq 2K+1}\sum_{i=1}^m||\partial_i F^\infty(x,\textbf{u})||_{2}||\textbf{u}_n-\textbf{v}_n-\textbf{w}_n||_{2}\\
&+\sup_{||\textbf{u}||_{H^\alpha}\leq 2K+1}\sum_{i=1}^m||\partial_i F^\infty(x,\textbf{u})||_{p}||\textbf{u}_n-\textbf{v}_n-\textbf{w}_n||_{p^*}\\
\leq&C_K(||\textbf{u}_n-\textbf{v}_n-\textbf{w}_n||_{2}+||\textbf{u}_n-\textbf{v}_n-\textbf{w}_n||_{p^*}),
\end{align*}
where $p^*$ is the conjugate of $p$.

Then, for any given $\epsilon>0$, there exists $N>0$ such that
\begin{align*}
J^\infty(\textbf{u}_n)-J^\infty(\textbf{v}_n)-J^\infty(\textbf{w}_n)\geq -\epsilon,\ \forall n\geq N.
\end{align*}
We have that
\begin{align*}
I^\infty_{c_1,\ldots,c_m}&=J^\infty(\textbf{u}_n)-\epsilon\\
&\geq J^\infty(\textbf{v}_n)+J^\infty(\textbf{w}_n)-2\epsilon\\
&\geq I^\infty_{a_{n,1},\ldots,a_{n,m}}+I^\infty_{b_{n,1},\ldots,b_{n,m}}-2\epsilon
%&\geq I^\infty_{\sqrt{c^2_{1}-a^2_{1}},\ldots,\sqrt{c^2_{m}-a^2_{m}}}+I^\infty_{a_{n,1},\ldots,a_{n,m}}-\epsilon.
\end{align*}
for $n$ large enough.

Since the map $(c_1,\ldots,c_m)\mapsto I^\infty_{c_1,\ldots,c_m}$ is continuous, let $n$ go to infinity. Thanks to the  arbitrariness of $\epsilon$, we have
\begin{align*}
I^\infty_{c_1,\ldots,c_m}\geq I^\infty_{c_{1}-a_{1},\ldots,c_{m}-a_{m}}+I^\infty_{a_{1},\ldots,a_{m}}
\end{align*}
which is a contraction with Lemma \ref{le2.3} .

Thus, we conclude that the compactness occurs.  Therefore, there exists a sequence $\{y_n\}_{n\in \mathbb{N}^+}\subset\mathbb{R}^N$ such that for any given $\epsilon>0$, there exists $r(\epsilon)>0$ with
\begin{align*}
\sum_{i=1}^m\int_{B(y_n,r(\epsilon))}u^2_{n,i}dx\geq c-\epsilon.
\end{align*}
For each $n\in\mathbb{N}^+$, we can choose $z_n=(z_{n,1},\ldots,z_{n,m})$ with $z_{n,i}\in\{NL,N=0,1,\ldots\}(1\leq i\leq m)$ such that $y_n-z_n\in[0,L]^N$, where $L\in\mathbb{R}^N$ is such that $F^\infty(x+L,\textbf{u})=F^\infty(x,\textbf{u})$ for all $x\in\mathbb{R}^N$ and $\textbf{u}\in\mathbb{R}^m$. Now, define $\textbf{v}_n=\textbf{u}_n(x+z_n)$. Note that $||\textbf{u}_n||_{H^\alpha}=||\textbf{v}_n||_{H^\alpha}$ which is bounded. Therefore, there exists a subsequence $\textbf{v}_n$ which is denoted by itself and there exists $\textbf{v}$ satisfying $\textbf{v}_n\rightharpoonup \textbf{v}$ in $H^\alpha$.

On the one hand,  there holds
\begin{align*}
\sum_{i=1}^m||v_i||^2_{2}&\geq \sum_{i=1}^m\int_{B(0,r+\sqrt{N}L)}|v_i|^2dx\\
&\geq \lim_{n\rightarrow\infty}\sum_{i=1}^m\int_{B(0,r+\sqrt{N}L)}|v_{n,i}|^2dx\\
&=\lim_{n\rightarrow\infty}\sum_{i=1}^m\int_{B(z_n,r+\sqrt{N}L)}|u_{n,i}|^2dx\\
&\geq \lim_{n\rightarrow\infty}\sum_{i=1}^m\int_{B(y_n,r)}|u_{n,i}|^2dx\\
&\geq c-\epsilon.
\end{align*}
Due to the  arbitrariness of $\epsilon$, we have
\begin{align*}
\sum_{i=1}^m||v_i||^2_{2}\geq c.
\end{align*}

On the other hand, $||v_i||^2_{2}\leq\liminf_{n\rightarrow\infty}||v_{n,i}||^2_{2}\leq c_i$ for $i=1,\ldots,m.$ Therefore, there holds
\begin{align*}
||v_i||^2_{2}=c_i
\end{align*}
and
\begin{align*}
\lim_{n\rightarrow\infty}||v_i-v_{n,i}||_{2}=0,\ 1\leq i\leq m.
\end{align*}
Applying  an interpolation theorem, we have
\begin{align*}
v_{n,i}\rightarrow v_i\ \textrm{in}\ L^p,\ p\in[2,2^*_\alpha).
\end{align*}
Then, we have that:
\begin{align}\label{3.1}
\int F^\infty(x,\textbf{v}_n)dx\rightarrow\int F^\infty(x,\textbf{v})dx.
\end{align}

Due to the periodicity of $F^\infty(\cdot,\textbf{u})$, we get
\begin{align*}
J^\infty(\textbf{u}_n)=J^\infty(\textbf{v}_n)\rightarrow I^\infty_{c_1,\ldots,c_m}\ \textrm{as}\ n\rightarrow\infty.
\end{align*}
By the weak lower semi-continuity of the norm in $L^2$ and \eqref{3.1}, it follows that
\begin{align*}
J^\infty(\textbf{v})=I^\infty_{c_1,\ldots,c_m}
\end{align*}
which implies
\begin{align*}
v_{n,i}\rightarrow v_i\ \textrm{in}\ H^\alpha,\ 1\leq i\leq m.
\end{align*}
 This completes the proof of Theorem \ref{th1.2}.

\section{Proof of Theorem \ref{th1.1}}

In this section, we will prove Theorem \ref{1.1}. Assume that $\textbf{u}_n$ is a minimizing sequence of the constrained variational problems of $I_{c_1,\ldots,c_m}$. First, we will rule out vanishing by contradiction. We suppose that $\textbf{u}_n$  is vanishing. That is,  for all $r<\infty$, it follows that
$$
\lim_{n\rightarrow\infty}\sup_{y\in \mathbb{R}^N}\sum_{i=1}^m\int_{|x-y|\leq r}|u_{n,i}|^{2}dx=0.
$$
From Lemma 2.4 in \cite{Guo}, one can see that $\lim_{n\rightarrow\infty}\sum_{i=1}^m||u_{n,i}||_p=0$ for $p\in(2,2^*_\alpha)$.

From the assumptions $(H_1)$ and $(H_4)$, it follows that for any given $\delta>0$, there exists $R_\delta>0$ such that
\begin{align*}
F(x,\textbf{u}_n)\leq\sum_{i=1}^m\Big[\delta(|u_{n,i}|^2+|u_{n,i}|^{l_{2,i}+2})+A'(|u_{n,i}|^{\beta'_i+2}+|u_i|^{l_{3,i}+2})\Big],\ \forall |x|\geq R_\delta.
\end{align*}
Therefore,
\begin{align}\label{4.1}
\int_{|x|\geq R_\delta} F(x,\textbf{u}_n)dx\leq\sum_{i=1}^m\Big[\delta(||u_{n,i}||_2^2+||u_{n,i}||^{l_{2,i}+2}_{l_{2,i}+2})+A'(||u_{n,i}||^{\beta'_i+2}_{\beta'_i+2}+||u_{n,i}||^{l_{3,i}+2}_{l_{3,i}+2})\Big].
\end{align}
This implies that
\begin{align*}
\limsup_{n\rightarrow\infty}\int_{|x|\geq R_\delta} F(x,\textbf{u}_n)dx\leq\delta c.
\end{align*}
Note that
\begin{align*}
&\int_{|x|\leq R_\delta} F(x,\textbf{u}_n)dx\\
\leq &A\sum_{i=1}^m\int_{|x|\leq R_\delta}\Big(|u_{i,n}|^2+|u_{i,n}|^{l_{3,i}+2}\Big)dx\\
\leq &A\sum_{i=1}^m(||u_{n,i}||^{l_{3,i}+2}_{l_{3,i}+2}|R_\delta|^{\frac{l_{3,i}}{l_{3,i}+2}}+||u_{n,i}||^{l_{3,i}+2}_{l_{3,i}+2})\rightarrow 0, \ \textrm{as}\ n\rightarrow\infty.
\end{align*}

From \eqref{4.1} and the  arbitrariness of $\delta$, we have
\begin{align*}
\lim_{n\rightarrow\infty}\int F(x,\textbf{u}_n)dx=0.
\end{align*}
Meanwhile, we notice that
\begin{align*}
\lim_{n\rightarrow\infty}J(\textbf{u})=I_{c_1\ldots,c_m}<0.
\end{align*}
Which is a contradiction.

Next, we will ruled out the dichotomy by contradiction.  Suppose that there exists constants $a_i\in(0,c_i)$ and two bounded vectorial sequences $\textbf{v}_{n}=(v_{n,1},\ldots,v_{n,m}),\textbf{w}_{n}=(w_{n,1},\ldots,w_{n,m})\subset H^{\alpha}$ such that for $1\leq i,j\leq m$, we have
$$
\mbox{supp}\ v_{n,i}\bigcap \mbox{supp}\ w_{n,j}=\emptyset,
$$
$$
|v_{n,i}|+|w_{n,i}|\leq |u_{n,i}|.
$$
\begin{align}\label{4.2}
||v_{n,i}||^{2}_{L^{2}}:=a_{n,i}\rightarrow a_i,\quad ||w_{n,i}||^{2}_{L^{2}}:=b_{n,i}\rightarrow (c_i-a_i),\  as\ n\rightarrow\infty.
\end{align}
$$
||u_{n,i}-v_{n,i}-w_{n,i}||_{L^{p}}\rightarrow 0,\ \textrm{for}\ 2\leq p<2^*_\alpha,
$$
where 
	\begin{align*}
	2^*_\alpha=\left\{
	\begin{array}{ll}
	\infty, & N\leq 2\alpha \\
	\frac{2N}{N-2\alpha}, & N> 2\alpha.
	\end{array}
	\right.
	\end{align*} 
\begin{eqnarray*}
\liminf_{n\rightarrow\infty}\{\langle(-\triangle)^{\alpha}u_{n,i},u_{n,i}\rangle-\langle(-\triangle)^{\alpha}v_{n,i},v_{n,i}\rangle-\langle(-\triangle)^{\alpha}w_{n,i},w_{n,i}\rangle\}\geq 0.
\end{eqnarray*}

Then, for any given $\epsilon>0$, choosing $\delta\in(0,\epsilon)$, there exists $n_2>0$ such that
\begin{align*}
&J(\textbf{u}_n)-J(\textbf{v}_n)-J^\infty(\textbf{w}_n)\\
=&\frac{1}{2}\sum_{i=1}^m\Big(\langle(-\triangle)^{\alpha}u_{n,i},u_{n,i}\rangle-\langle(-\triangle)^{\alpha}v_{n,i},v_{n,i}\rangle-\langle(-\triangle)^{\alpha}w_{n,i},w_{n,i}\rangle\Big)\\
&-\int(F(x,\textbf{u}_n)-F(x,\textbf{v}_n)-F(x,\textbf{w}_n))dx+\int(F^\infty(x,\textbf{w}_n)-F(x,\textbf{w}_n))dx\\
\geq&-\delta-\int (F(x,\textbf{u}_n)-F(x,\textbf{v}_n+\textbf{w}_n))dx+\int(F^\infty(x,\textbf{w}_n)-F(x,\textbf{w}_n))dx\\
\geq&-\delta-\int (F(x,\textbf{u}_n)-F(x,\textbf{v}_n+\textbf{w}_n))dx+\int_{|x-y_n|\geq R_n}(F^\infty(x,\textbf{w}_n)-F(x,\textbf{w}_n))dx,\ \forall \ n\geq n_2,
\end{align*}
where $\{y_n\}_{n\in\mathbb{N}^+}$ is a bounded sequence in $\mathbb{R}^N$, $R_n$ satisfies $\lim_{n\rightarrow\infty}R_n=+\infty$.

With similar arguments as in the proof of theorem \ref{th1.2}, there exists $n_3>0$ such that
\begin{align*}
-\delta-\int(F(x,\textbf{u}_n)-F(x,\textbf{v}_n+\textbf{w}_n))dx\geq -\epsilon, \ \forall\ n\geq n_3.
\end{align*}

Therefore, it follows that
\begin{align*}
J(\textbf{u}_n)-J(\textbf{v}_n)-J^\infty(\textbf{w}_n)\geq -\epsilon+\int_{|x-y_n|\geq R_n}(F^\infty(x,\textbf{w}_n)-F(x,\textbf{w}_n))dx,\ \forall\ n\geq\max\{n_2,n_3\}.
\end{align*}

From the assumption $(H_4)$, for any given $\eta>0$, there exists $R>0$ such that
\begin{align*}
|F^\infty(x,\textbf{u})-F(x,\textbf{u})|\leq \eta\sum_{i=1}^m(|u_i|^2+|u_i|^{l_{2,i}+2})
\end{align*}
for all $\textbf{u}$ and $|x|\geq R.$

Thanks to the boundness of $\{y_n\}$ and $\lim_{n\rightarrow\infty}R_n=\infty$, there exists $n_4>0$ such that
\begin{align*}
\{x:|x-y_n|\geq R_n\}\subset\{x:|x|\geq R\},\ \forall\ n\geq n_4.
\end{align*}
Then, we have
\begin{align*}
\int_{|x-y_n|\geq R_n}|F^\infty(x,\textbf{w}_n)-F(x,\textbf{w}_n)|dx\leq \eta\sum_{i=1}^m(\|w_{n,i}\|_2^2+\|w_{n,i}\|_{l_{2,i}+2}^{l_{2,i}+2}),\ \forall\ n\geq n_4.
\end{align*}
Combining that and the boundedness of $\{\textbf{w}_n\}$ in $H^\alpha$, there holds
\begin{align*}
\lim_{n\rightarrow\infty}\int_{|x-y_n|\geq R_n}|F^\infty(x,\textbf{w}_n)-F(x,\textbf{w}_n)|dx=0.
\end{align*}

Therefore, it obtains that
\begin{align*}
I_{c_1,\ldots,c_m}&\geq \lim_{n\rightarrow\infty}J(\textbf{u}_n)\\
&\geq \liminf_{n\rightarrow\infty}\Big(J(\textbf{v}_n)+J^\infty(\textbf{w}_n)\Big)-\epsilon\\
&\geq \liminf_{n\rightarrow\infty}\Big(I_{a_{n,1},\ldots,a_{n,m}}+I^\infty_{b_{n,1},\ldots,b_{n,m}}\Big)-\epsilon.
\end{align*}

From \eqref{4.2} and Lemma \ref{le2.1}, $I_{c_1,\ldots,c_m}\geq I_{a_1,\ldots,a_m}+I^\infty_{c_1-a_1,\ldots,c_m-a_m}-\epsilon$.  Thanks to the  arbitrariness of $\epsilon$, we get
\begin{align*}
I_{c_1,\ldots,c_m}\geq I_{a_1,\ldots,a_m}+I^\infty_{c_1-a_1,\ldots,c_m-a_m}.
\end{align*}
It is a contradiction which implies that the dichotomy cannot occur.

Thus, we can conclude that the compactness occurs. Hence, there exists a sequence $\{y_n\}_{n\in\mathbb{N}^+}\subset \mathbb{R}^N$ such that for all $\epsilon>0$, there exists $R(\epsilon)>0$ with
\begin{align}\label{4.3}
\sum_{i=1}^m\int_{B(y_n,R(\epsilon))}u^2_{n,i}dx\geq c-\epsilon.
\end{align}

Now we will prove the sequence $\{y_n\}_{n\in\mathbb{N}^+}$ is bounded. By contradiction, we suppose that $|y_n|\rightarrow\infty$ by passing to a subsequence which is denoted by itself. We can choose $z_n=(z_{n,1},\ldots,z_{n,m})$ satisfying $z_{n,i}\in\{NL,N=0,1,\ldots\}$, $\lim_{n\rightarrow\infty}|z_{n,i}|=\infty$, $i=1,\ldots,m$ such that $y_n-z_n\in [0,L]^N$.

Define $\textbf{v}_n(x)=\textbf{u}_n(x+z_n):=(v_{n,1},\ldots,v_{n,m})$. With a similar argument in the proof of Theorem \ref{th1.2}, there exists $\textbf{v}=(v_1,\ldots,v_m)\in H^\alpha$ such that $v_{n,i}\rightharpoonup v_i$ in $H^\alpha$ and
\begin{align*}
\lim_{n\rightarrow\infty}||v_{n,i}-v_i||_{L^p}=0\ \textrm{for}\ 2\leq p< 2^*_\alpha,\ i=1,\ldots,m.
\end{align*}
%\begin{align*}
%J^\infty(\textbf{v}_n)=J^\infty(\textbf{u}_n).
%\end{align*}

On the other hand, note that
\begin{align*}
J(\textbf{u}_n)-J^\infty(\textbf{u}_n)&=\int(F^\infty(x,\textbf{u}_n)-F(x,\textbf{u}_n))dx\\
&=\int(F^\infty(x,\textbf{u}_n)-F(x-z_n,\textbf{v}_n))dx.
\end{align*}
From the assumption $(H_4)$, for any given $\delta>0$, there exists $R(\delta)>0$ such that
\begin{align*}
&\Big|\int_{|x-z_n|\geq R(\delta)}\big(F^\infty(x,\textbf{v}_n)-F(x-z_n,\textbf{v}_n)\big)dx\Big|\\
=&\Big|\int_{|x-z_n|\geq R(\delta)}\big(F^\infty(x-z_n,\textbf{v}_n)-F(x-z_n,\textbf{v}_n)\big)dx \Big|\\
\leq&\delta\sum_{i=1}^m\int_{|x-z_n|\geq R(\delta)}(|v_{n,i}|^2+|v_{n,i}|^{l_{2,i}+2})dx\\
\leq&C\delta\sum_{i=1}^m\Big(||v_{n,i}||^2_{H^1}+||v_{n,i}||^{l_{2,i}+2}_{l_{2,i}+2}\Big)\\
\leq&C\delta.
\end{align*}
The last inequality is from the Sobolev embedding and the boundedness of $\textbf{v}_n$ in $H^\alpha$.

Since $\lim_{n\rightarrow\infty}|z_n|=\infty$, there exists $n_{R(\delta)}>0$ such that
\begin{align*}
&\Big|\int_{|x-z_n|\leq R(\delta)}F^\infty(x,\textbf{v}_n)-F(x-z_n,\textbf{v}_n)dx\Big|\\
\leq&\Big|\int_{|x|\geq \frac{1}{2}|z_n|}F^\infty(x,\textbf{v}_n)-F(x-z_n,\textbf{v}_n)dx\Big|\\
\leq& K\sum_{i=1}^m\int_{|x|\geq \frac{1}{2}|z_n|}\Big(|v_{n,i}|^2+|v_{n,i}|^{l_{2,i}+2}\Big)dx\\
\leq& K\sum_{i=1}^m\Big(\int_{|x|\geq \frac{1}{2}|z_n|}\big(|v_i|^2+|v_i|^{l_{2,i}+2}\big)dx+\int_{|x|\geq \frac{1}{2}|z_n|}\big(|v_{n,i}-v_i|^2+|v_{n,i}-v_i|^{l_{2,i}+2}\big)dx\Big)\\
\leq& K\sum_{i=1}^m\Big(\int_{|x|\geq \frac{1}{2}|z_n|}\big(|v_i|^2+|v_i|^{l_{2,i}+2}\big)dx+\int\big(|v_{n,i}-v_i|^2+|v_{n,i}-v_i|^{l_{2,i}+2}\big)dx\Big)
\end{align*}
for $n\geq n_{R(\delta)}.$

It obtains that
\begin{align*}
\lim_{n\rightarrow\infty}\Big|\int_{|x-z_n|\leq R}F^\infty(x,\textbf{v}_n)-F(x-z_n,\textbf{v}_n)dx\Big|=0.
\end{align*}

Therefore,  we conclude that
\begin{align*}
\liminf_{n\rightarrow\infty}(J(\textbf{u}_n)-J^\infty(\textbf{u}_n))\geq -C\delta.
\end{align*}
Due to the arbitrariness of $\delta,$ it obtains that
\begin{align*}
I_{c_1,\ldots,c_m}=\lim_{n\rightarrow\infty}J(\textbf{u}_n)\geq\liminf_{n\rightarrow\infty}J^\infty(\textbf{u}_n)\geq I^\infty_{c_1,\ldots,c_m}.
\end{align*}
It is a contradiction with Lemma \ref{le2.4}. Hence,  $\{y_n\}$ is bounded.

Let $\rho=\sup_{n\in\mathbb{N}}|y_n|$. From \eqref{4.3}, there holds
\begin{align*}
\sum_{i=1}^m\int_{B(0,R(\epsilon)+\rho)}|u_{n,i}|^2 dx\geq \sum_{i=1}^m\int_{B(y_n,R(\epsilon))}|u_{n,i}|^2dx\geq c-\epsilon.
\end{align*}
Therefore, we have
\begin{align*}
\sum_{i=1}^m\int|u_{i}|^2dx\geq \sum_{i=1}^m\int_{B(0,R(\epsilon)+\rho)}|u_{i}|^2dx=\lim_{n\rightarrow\infty}\sum_{i=1}^m\int_{B(0,R(\epsilon)+\rho)}|u_{n,i}|^2dx\geq c-\epsilon.
\end{align*}
Due to the arbitrariness of $\epsilon$, it obtains that
\begin{align*}
\sum_{i=1}^m\int |u_i|^2dx\geq c.
\end{align*}

From the weakly lower semicontinuity of the $L^2$ norm, one can see
\begin{align*}
\int |u_i|^2dx\leq c^2_i,\ \forall 1\leq i\leq m.
\end{align*}
Therefore, $\textbf{u}\in S_c$ and $\lim_{n\rightarrow\infty}||u_{n,i}-u_i||_2=0$, $1\leq i\leq m$. From the boundedness of $||\textbf{u}_n||_{H^1}$ and Gagliardo-Nirenberg inequality, it obtains that $\lim_{n\rightarrow\infty}||u_{n,i}-u_i||_p=0$ for $p\in [2,2^*_\alpha)$, $1\leq i\leq m$. It implies that $\lim_{n\rightarrow\infty}\int F(x,\textbf{u}_n)dx=\int F(x,\textbf{u})dx$ and $J(\textbf{u})=I_{c_1,\ldots,c_m}$.

\end {document}